
\input amstex
\documentstyle{amsppt}
\magnification=1200
\parindent 20 pt
\nopagenumbers
\NoBlackBoxes

\define \a{\alpha}
\define \be{\beta}

\define \dl{\delta}

\define \g{\gamma}
\define \G{\Gamma}

\define \lm{\lambda}

\define \om{\omega}

\define \th{\theta}

\define \ve{\varepsilon}
\define \vp{\varphi}

\define \z{\zeta}

\define \df{\dsize\frac}

\define \fc{\frac}
\define \iy{\infty}

\define \sq{\sqrt}
\define \tm{\times}

\define \ri{\rightarrow}

\define \pf{\demo{Proof}}
\define \edm{\enddemo}
\define \ep{\endproclaim}

\define \sk{\smallskip}
\define \mk{\medskip}
\define \bk{\bigskip}

\define \1{^{-1}}
\define \2{^{-2}}

\define \hf{\hat{f}}

\define \BR{\Bbb R}

\define \tf{\tilde{f}}

\define \Lip{\operatorname{Lip}}

\topmatter
\centerline{\bf SOME CONDITIONS FOR EXISTENCE AND INTEGRABILITY}
\centerline{\bf OF THE FOURIER TRANSFORM }
\vskip .20in
\centerline{by}
\vskip .20in
\centerline{\bf E. Liflyand}
\centerline{\it Dept. of Mathematics \& Computer Science}
\centerline{\it Bar-Ilan University}
\centerline{\it Ramat-Gan 52900}
\centerline{ISRAEL}

\baselineskip 20pt
\vskip .30in
\abstract{The Fourier transform is naturally defined for integrable functrions.
Otherwise, it should be stipulated in which sense the Fourier transform is
understood.
We consider some class of radial and, generally saying, nonintegrable
functions.
The Fourier transform is calculated as an improper integral and the limit
coincides with the Fourier transform in the distributional sense.
The inverse Fourier formula is proved as well.
Given are some applications of the result obtained.}  \endabstract
\endtopmatter
\baselineskip 20pt
\document
\vskip 3cm
\centerline{\bf Preprint\qquad BIMACS--9502}
\vskip 4cm
\centerline{\bf Bar-Ilan University,\qquad 1995}
\vfill \eject
\subheading{1. INTRODUCTION}

{\bf 1.1.} Let $\BR^n$ be a real Euclidean $n$-dimensional space with elements
$x=(x_1,\dots, x_n).$
If a function $f(x)$ is integrable, in the Lebesgue sense, over all $\BR^n$
there is no problem to define its Fourier transform
$$\hf(x) = \int_{\BR^n} f(u)\ e^{-ix\cdot u} du$$
where $x\cdot u=x_1 u_1 +\dots + x_n u_n$ is the scalar product of
$u,x\in\BR^n.$ But, even in this case, special additional conditions are
needed for the inverse formula to be true:
$$f(u) = (2\pi)^{-n} \int_{\BR^n} \hf (x)\ e^{iu\cdot x} dx.$$
For instance, sometimes a summability method should be applied to impart
a certain sense to the last integral.
And if $f$ is not integrable both formulas need a special investigation.

Thus, the main purpose of our work is to give certain conditions under which
both formulas hold.
We will formulate them just after introducing notations needed.
Then a separate section will be devoted to auxiliary lemmas.
The next section will present the proof of the theorem.
After that, several sections will be devoted to applications.
These applications were announced in \cite{BL1}, but the proofs were given only
in the author's Ph.D. thesis (see also [L2] where conditions are more
restrictive than in this paper).
The first application gives a connection between many-dimensional and
one-dimensional Fourier transforms.
This allows to reduce some multidimensional problems to easier or known
one-dimensional ones.
Then we give an asymptotic formula for the Fourier transform of a radial
function with convexity type conditions.

{\bf 1.2.} We consider radial functions $f(x) = f_0(|x|),$ where $|x|=
(x\cdot x)^{\fc{1}{2}}$ is not necessarily integrable over all $\BR^n.$
Some notations are needed to describe precisely the class of functions to be
studied. For the sake of convenience and completeness, some other notations,
which will be used in proofs, are given here as well.

Let $0<\a\leq {n-1\over 2}$ and $\a^*$ be the greatest integer smaller than
$\a.$ Denote by
$$W_\a (f_0 ; t) = \fc{1}{\G(\a)} \int^\iy_t f_0(r)(r-t)^{\a-1} dr$$
the Weyl integral of fractional order where $r$ and $t$ are real numbers.
When\linebreak
$0<\a<1$
$$f_0^{(\a)} (t) =\fc{d}{dt}\ W_{1-\a}\ (f_0 ; t)$$
is the Weyl derivative of fractional order.
When $\a=p+\g$ with $p=1,2,\dots\ $, and $0<\g<1,$ then
$$f_0^{(\a)} (t) = \fc{d^p}{dt^p}\ f_0^{(\g)}\ (t).$$
Let
$$R_\a (f_0 ; t) = \fc{1}{\G(\a)} \int^t_0 f_0 (r)\ (t-r)^{\a-1} dr$$
be the Riemann-Liouville fractional integral.
All these notions may be found for example in \cite{BE2}, Ch.13 (see
also [Co], [SKM], [Tr]).

Let $C[a,b]$ and $C^p [a,b]$ be the classes of continuous functions and
having $p$ continuous derivatives, respectively.

Let us introduce two Bessel-type functions:
$$\align Q_\a (t) &= \int^1_0 (1-s)^{\a-1} s^{\fc{n}{2}}\
J_{\fc{n}{2}-1} (ts) ds \\ \quad\\
&=\G(\a) t^{-{n\over 2}-\a} R_\a(s^{{n\over 2}}J_{{n\over 2}-1}(s);t)
\endalign$$
and
$$q_\a (t) = \int^1_0 (1-s)^{\a-1} s^{\fc{n}{2}-1}\ J_{\fc{n}{2}} (ts) ds$$
where
$J_\mu$ is the Bessel function of first type and order $\mu.$

In what follows $F_\a (t) = t^{\fc{n-1}{2}}\ f_0^{(\a)}\ (t).$
When $\a={n-1\over 2},$ for brevity, we will write simply $Q,q,F.$

We define $\varphi\in S$ when $\varphi\in C^\infty$ and $\varphi(x)$
and its derivatives, which we are allowed to multiply by any polynomial
of $|x|,$ tend to $0$ uniformly as $|x|\to\infty.$ Continuous linear
forms on $S$ are called tempered distributions.

We shall denote absolute constants, that is independent of substantial
parameters, by the letter $C.$ If a dependence on some parameters is
essential, they will be indicated as subscripts. Whenever no confusion
can result, we use the same letter for different constants in
different places.

{\bf 1.3.} Let us consider a class of radial functions satisfying the
following conditions:
$$\align &f_0,..., f^{(\a^*)}\quad\text{are locally absolutely
continuous on}\quad (0,\iy);\tag1\\
&\lim_{t\ri\iy} t^p
f_0^{(p)} (t) = 0 \quad\text{for}\quad p=0,1,\dots,\a^*; \tag2\\
&\lim_{t\ri\iy} F_\a (t) = 0; \tag3\\
&F_\a \quad \text{is a function of bounded variation
on}\quad (0,\iy).\tag4\endalign$$
The total variation of the function $F_{\a}$ will be denoted by $V_{F_{\a}}.$
\sk
\proclaim{Theorem 1}\ Let a function $f$ be radial and satisfy conditions
{\rm (1) - (4).} Then for $|x|>0$
$$f(x) = \lim_{A\ri\iy} (2\pi)^{-n}\int_{|u|\leq A}
\left(1-\fc{|u|^2}{A^2}\right)^{\fc{n-1}{2}-\a} \hf(u)\ e^{ix\cdot u} du\tag5$$
where
$$\hf(u) = \fc{(2\pi)^{\fc{n}{2}}(-1)^{\a^*+1}}{\G(\a)}\ |u|^{1-\fc{n}{2}}
\int^\iy_0 F_\a (t)\ t^{\a+\fc{1}{2}} Q_\a (|u| t)\ dt\tag6$$
is continuous, tends to zero as $|x|\ri\iy$ and coincides with the Fourier
transform $\tf,$ of the function $f,$ in the distributional sense.
Both integrals converge uniformly for $|x|\geq r_0>0.$\ep

We have to mention that the case $\a={n-1\over 2}$ was considered earlier in
\cite{BL1}, \cite{BL2}.
Some ideas of the proof in \cite{BL2} were used here and also in [L2] where
such results are established under superfluous conditions at zero.

{\bf 1.4.} Compare Theorem 1 with some earlier results.
S. Bochner in \cite{Bo}, \S 44, considers not only radial functions, but more
restrictive conditions are claimed for the radial part of a function (spherical
average of a function).
Analogously, the radiality allows less restrictive smoothness conditions than
those in V.A. Ilyin and S. A. Alimov theorems for general spectral expansions
(see \cite{IA}). In M. L. Goldman's paper \cite{G} radial functions are
considered, with ``worse'' conditions at infinity and the monotonicity of a
given function and its derivatives. Very simple formula somehow similar to
(6) may be found in [SKM], Ch.5, Lemma 25.1', but the authors did not care
much for sharp assumptions.
\mk

\subheading{2. AUXILIARY LEMMAS}
\sk
\proclaim{Lemma 1} The following asymptotic relation holds
$$q_\a(r) = \G(\a)\ r^{-\a} J_{\fc{n}{2} + \a} (r) + \z_{\a,n} r^{-\fc{n}{2}}
+O \left(r^{-\a-\fc{3}{2}}\right)$$
as $r\ri\iy,$ and $\z_{\a,n}$ are some numbers. \ep
\sk
\demo{Proof}\ We have
$$\align q_\a(r) &= \sum\limits_{j=0}^M \int^1_0 (1-s)^{\a-1+j}
(1+s)^j s^{\fc{n}{2}+1} J_{\fc{n}{2}} (rs)\ ds\\
&+ \int^1_0 g_0(s)\ s^{\fc{n}{2}-1}\ J_{\fc{n}{2}} (rs)\,ds\endalign$$
where $M$ is such that $g_0(s)=(1-s)^{\a+M}(1+s)^{M+1}$ is smooth enough.
Evaluate first the last integral. We need the following properties of the
Bessel functions (see e.g., [BE2], \S 7.2.8(50),(51); \S 7.13.1(3);
\S 7.12(8)):
$$\align
&\fc{d}{dt}\ \left[ t^{\pm\nu} J_\nu (t)\right] = \pm t^{\pm\nu}\ J_{\nu\mp 1}\
(t);\tag7\\
&J_\nu (t) = \sqrt{\fc{2}{\pi t}} \cos\left(t-\fc{\pi\nu}{2}-\fc{\pi}{4}
\right)+ \sqrt{\fc{2}{\pi}}\ \fc{1-4\nu^2}{8}\ t^{-\fc{3}{2}} \sin
\left(t-\fc{\pi\nu}{2}-\fc{\pi}{4}\right)\tag8\\ \quad\\
&\qquad\ + O\left(t^{-\fc{5}{2}}\right)\qquad \text{as}\quad t\ri\iy;\\
\quad\\&J_\nu = O (t^\nu)\qquad \text{for small}\quad t.\tag9\endalign$$
Let us integrate by parts, using (7) as follows:
$$\align
\int^1_0  g_0(s)\ s^{\fc{n}{2}-1}\ J_{\fc{n}{2}} (rs)\ ds&=
\int^1_0 [g_0(s)\ s^{n-2}]\ [s^{-\fc{n}{2}+1}\ J_{\fc{n}{2}} (rs)]\ ds\\
=-{1\over r} g_0(s)\ s^{\fc{n}{2}-1}\ J_{\fc{n}{2}-1} (rs)\biggm|_0^1
&+{1\over r}\int^1_0 g_1(s)\ s^{\fc{n}{2}-2}\ J_{\fc{n}{2}-1} (rs)\
ds\endalign$$ where $g_1$ is also smooth enough. For the case $n=2,$
unlike for that for $n>2$ the factor $s^{{n\over 2}-2}$ does not
appear in the last integral. We can continue this procedure if needed.
After $[{n\over 2}]$ steps we get
$$\z_{\a,n} r^{-{n\over 2}}+r^{-{n\over 2}}\int_0^1 g_{[{n\over 2}]}
(s) J_0(rs)\,ds=\z_{\a,n} r^{-{n\over 2}}+O(r^{-{n+2\over 2}})$$
when $n$ is even, and the integral
$$r^{-{n-1\over 2}}\int_0^1 g_{[{n\over 2}]}(s) s^{{1\over 2}}
J_{{1\over 2}}(rs)\,ds=\z_{\a,n} r^{-{n\over 2}}+O(r^{-{n+2\over 2}})$$
when $n$ is odd.

To estimate the sum we need the following lemma.
\sk
\proclaim{Lemma 2}\ For $r\geq 1,\ \be>-\df{1}{2},\ \mu > -1$ and each
positive integer $p$
$$\int^1_0 (1-s)^\mu\ s^{\be+1} J_\be (rs)\ ds = \sum^p_{j=1} \a^\mu_j
r^{-(\mu+j)}\ J_{\be+\mu+j} (r) + O\left(r^{-\mu-p-\fc{3}{2}}\right),$$
where $\a^\mu_j$ are some numbers depending only on $p$ and $\mu,$ and
$\a^\mu_1 = \G(\mu+1),\ \a^\mu_2 = \mu\G(\mu+2).$\ep
\sk
\demo{Proof}\ We have
$$\align (1-s)^\mu &= 2^{-\mu} (1-s^2)^\mu + \left[(1-s)^\mu - 2^{-\mu}
(1-s^2)^\mu\right]\\ \quad\\
 &= 2^{-\mu} (1-s^2)^\mu + 2^{-\mu} (1-s)^{\mu+1}\ \fc{2^\mu-(1+s)^\mu}
{1-s}.\endalign$$
Continuing this process of chipping off the binomials $(1-s^2)^{\mu+j-1}$
for $j=1,\dots, p,$ use the formula (see e.g., \cite{SW}, Ch.4, Lemma 4.13):
$$J_{\be+\mu+1} (r) = \fc{r^{\mu+j}}{2^{\mu+j-1}\G(\mu+j)}\ \int^1_0 J_\be
(rs)\ s^{\be+1} (1-s^2)^{\mu+j-1} ds.$$

The remainder term is estimated as above, for the second integral for
$q_\a(r),$ by integrating  by  parts $\mu^* + p + 1$ times. Estimates
are better in this case since $s$ here is in rather high power.
The lemma is proved.\qquad \qed \edm

\sk

To finish the proof of Lemma 1 it remains to apply Lemma 2, with
$\be=\df{n}{2},\ \mu=\a-1+j,\ p=1,$ to the integrals in the sum
for $q_\a(r).$ The proof of Lemma 1 is complete.\qquad \qed \edm
\sk
\remark{Remark 2}\ Sometimes the estimate $q_\a(r) =
O\left(r^{-\a-\fc{1}{2}}\right)$ will be enough.\endremark
\sk
The following lemma is due to Trigub (see \cite{T4}, Lemma 2).
Since that edition is difficult of access, we give the proof here.
The lemma deals with the functions
$$i(\mu,\lm,r) = \int^1_0 t^\mu J_\lm (rt)\ dt$$
where $\mu+\lm> -1.$
\sk
\proclaim{Lemma 3}
\roster
\item "1)" $i(\mu,\lm,r) = \df{1}{r}\ J_{\lm+1}(r) + \df{\lm+1-\mu}{r}\
i(\mu-1,\lm + 1,r).$
\item "2)" The function $i(\mu,\lm,r)$ is $O(r^\lm)$ for small $r,$ and when
$r\ri\iy$ it behaves as $O\left(r^{-{3\over 2}}\right)$ or
$O\left(r^{-1-\mu}\right)$ for $\mu > {1\over 2}$ and $\mu \leq
{1\over 2},$ respectively.\endroster\ep
\sk
\demo{Proof}\ To prove 1), integrate by parts using (7).
We have
$$\align
i(\mu,\lm,r) &= \df{1}{r} \int^1_0 t^{\mu-\lm-1} d\left[ t^{\lm+1}
J_{\lm+1}(rt)\right]\\ \quad\\
&= \df{1}{r}\ J_{\lm+1}(r) + \df{\lm+1-\mu}{r} \int^1_0 t^{\mu -1}
J_{\lm+1}(rt)\ dt.\endalign$$
The first assertion in 2) immediately follows from (9) applied to $J_\lm$
and the usual estimate of the integral. Now let $r\geq 1.$
After a linear change of variables
$$i(\mu,\lm,r) = r^{-1-\mu} \int^r_0 t^\mu J_\lm (t)\ dt.$$
Decompose the integral. It is bounded when $t\in [0, 1].$ When $t\in [1, r]$
use (8). If $\mu > \fc{1}{2},$ then we obtain after integrating by parts
$$\align
& \int^r_1 t^\mu J_\lm (t)\ dt
=\int^r_1 t^\mu \left[ \sqrt{\fc{2}{\pi}}\ \fc{\cos\left(t-\fc{\pi\lm}{2} -
\fc{\pi}{4}\right)} {\sqrt{t}} + O\left(t^{-\fc{3}{2}}\right)\right]\
dt\\ \quad\\
&=\sqrt{\fc{2}{\pi}}\ r^{\mu-\fc{1}{2}} \sin \left( r-\fc{\pi\lm}{2} -
\fc{\pi}{4}\right) + O(1) + \int^r_1 O\left( t^{\mu-\fc{3}{2}}\right) dt
= O\left(r^{\mu-\fc{1}{2}}\right),\endalign$$
and
$$i(\mu,\lm,r) = O \left( r^{-1-\mu} r^{\mu-\fc{1}{2}}\right) = O
\left(r^{-\fc{3}{2}}\right).$$
If $\mu\leq\fc{1}{2},$ then the same computations show that the
integral $\int^r_1 t^\mu J_\lm(t) dt$ is bounded with respect to $r$
and
$$i(\mu,\lm,r) = O \left( r^{-1-\mu}\right).$$
The lemma is proved.\qquad \qed\edm

\bk

\subheading{3. PROOF OF THEOREM 1}
\sk
{\bf 3.1.} Let us start with proving (6).
Due to \cite{M}, p. 140, for each $\vp\in S$
$$\align
(\tf,\vp) &= \int_{\BR^n} \tf (x)\ \vp (x)\ dx\\
&=\lim_{A\ri\iy} \int_{\BR^n} \vp (x) \left[ (2\pi)^{\fc{n}{2}}
|x|^{1-\fc{n}{2}} \int^A_0 f_0 (t)\ t^{\fc{n}{2}} J_{\fc{n}{2}-1} (|x|
t)\ dt\right] dx.\endalign$$
This equality may be rewritten as follows:
$$(\tf,\vp) = (2\pi)^{\fc{n}{2}} \lim_{A\ri\iy} \int^\iy_0 \vp_0 (r)
\ r^{\fc{n}{2}} dr \int^A_0 f_0 (t)\ t^{\fc{n}{2}} J_{\fc{n}{2}-1} (rt)\ dt.$$
To prove (6) we will use the following procedure: integration by parts
so many times that the right-hand side of (6) is obtained, and then
justification of the passage to limit under the integral sign.
Thus, for $[\a]\ne 0$ integrate in $t$ by parts $[\a]$ times. We obtain
$$\align & \text{(10)} \qquad \qquad \qquad \qquad
\ \ \int^\iy_0 \vp_0 (r)\ r^{\fc{n}{2}} dr \int^A_0 f_0 (t)\ t^{\fc{n}{2}}
J_{\fc{n}{2}-1} (rt)\ dt \\ \quad\\
&\ \ = \int^\iy_0 \vp_0 (r)\ r^{\fc{n}{2}} \left\{ \fc{(-1)^{[\a]}}
{([\a]-1)!} \int^A_0 f_0^{([\a])} (t)\ t^{[\a]+\fc{n}{2}} dt \int^1_0
(1-s)^{[\a]-1} s^{\fc{n}{2}} J_{\fc{n}{2}-1} (rts)\ ds\right.\\ \quad\\
&\ \ \left. +\sum^{[\a]-1}_{p=0} \fc{(-1)^p}{p!}\ f_0^{(p)} (t)\
t^{p+\fc{n}{2}+1} \ \int^1_0 (1-s)^{p} s^{\fc{n}{2}} J_{\fc{n}{2}-1} (rts)\ ds
 \biggm|^A_0 \right\} dr.\endalign$$
To show that integrated terms vanish when $t=0$ apply (9)
to $J_{{n\over 2}-1}(rts).$
Let us demonstrate that $$\lim\limits_{t\to0+}f_0^{(p)}(t)t^{p+n}=0\qquad
\text{for}\qquad p=0,1,...,[\a]-1\qquad \text{and}\qquad \a\ge 1.$$
For $\a$ fractional we have $$f_0^{([\a])}(t)={1\over\Gamma
(\a-[\a])}\int\limits_t^\infty (s-t)^{\a-[\a]-1}f_0^{(\a)}(s)\,ds.$$
This formula due to Cossar [Co] may be found also in [Tr], Lemma 3.10.
It is easy to see now that $$f_0^{(p)}(t)={(-1)^{[\a]-p}\over\Gamma
(\a-p)}\int\limits_t^\infty (s-t)^{\a-p-1}f_0^{(\a)}(s)\,ds.$$
It is clear that for $\a$ integer the same formula holds. Integrate by parts
as follows $$\align \int\limits_t^\infty (s-t)^{\a-p-1}f_0^{(\a)}(s)\,ds
&=-{1\over\a-p} (s-t)^{\a-p}f_0^{(\a)}(s)\,ds \biggm|_t^\infty\\&+
{1\over\a-p}\int\limits_t^\infty (s-t)^{\a-p}f_0^{(\a+1)}(s)
\,ds\\&={1\over\a-p}\int\limits_t^\infty (s-t)^{\a-p}{1\over
s^{{n-1\over 2}}}\,dF_{\a}(s)\\&-{1\over\a-p}{n-1\over 2}\int\limits_
t^\infty (s-t)^{\a-p}{1\over s}f_0^{(\a)}(s)\,ds.\endalign$$
Here and below (3) is taken into account for integrated terms.
Split the last integral into two one of which is the same as that on the
left-hand side. Taking into account the corresponding coefficient we get
$$\align\biggl(1+{1\over\a-p}{n-1\over 2}\biggr)\int\limits_t^\infty (s-t)^
{\a-p-1}f_0^{(\a)}(s)\,ds&={1\over\a-p}\int\limits_t^\infty (s-t)^{\a-p}
{1\over s^{{n-1\over 2}}}\,dF_{\a}(s)\\&+{t\over\a-p}{n-1\over 2}
\int\limits_t^\infty (s-t)^{\a-p-1}{1\over s}f_0^{(\a)}(s)\,ds.\endalign$$
The last integral is equal to     $$\int\limits_t^\infty
(s-t)^{\a-p-1}{1\over s^{{n+1\over 2}}}F_{\a}(s)\,ds.$$
Integrate by parts again and obtain $$\align \int\limits_t^\infty
(s-t)^{\a-p-1}{1\over s^{{n+1\over 2}}}F_{\a}(s)\,ds&=-F_{\a}(s)
\int\limits_s^\infty (u-t)^{\a-p-1}{1\over u^{{n+1\over 2}}}\,du
\biggm|_t^\infty\\
&+\int\limits_t^\infty\,dF_{\a}(s)\int\limits_s^\infty (u-t)^{\a-p-1}
{1\over u^{{n+1\over 2}}}\,du\\ &= F_{\a}(t)
\int\limits_t^\infty (u-t)^{\a-p-1}{1\over u^{{n+1\over 2}}}\,du\\
&+\int\limits_t^\infty\,dF_{\a}(s)\int\limits_s^\infty (u-t)^{\a-p-1}
{1\over u^{{n+1\over 2}}}\,du.\endalign$$  Let us consider the integral
$$\int\limits_s^\infty (u-t)^{\a-p-1}{1\over u^{{n+1\over 2}}}\,du=
\int\limits_s^\infty u^{-{n+1\over 2}+\a-p-1}\,du.$$
It can be obviously estimated by $t^{-{n+1\over 2}+\a-p}$ and now
condition (4) completes the proof.

These calculations may seem somewhat superfluous. Indeed, for most cases
an easier way gives the result needed but for $p=0$ when $\a={n-1\over 2}$
one has to be more careful.

For $t=A$, taking into account that $A\ri\iy,$ it suffices,
in view of (2), to establish the uniform boundedness in $t$ of the integrals
$$B_p(t) = \int^\iy_0 \vp_0 (r)\ r^{\fc{n}{2}} dr t^{\fc{n}{2}+1}
\int^1_0 (1-s)^{p} s^{\fc{n}{2}} J_{\fc{n}{2}-1} (rts)\ ds.$$
Integrate by parts in the outer integral $m=\left[{n-1\over 2}\right]$ times,
using (7) so that the order of the Bessel function decreases.
Integrated terms vanish since $\vp_0\in S.$
We have (defining by $\psi$ here and below some function from $S$)
$$B_p(t) = \int^\iy_0 \psi (r)\ r^{\fc{n}{2}-m} t^{\fc{n}{2}-m+1} dr
\int^1_0 (1-s)^{p} s^{\fc{n}{2}-m} J_{\fc{n}{2}-m-1} (rts)\ ds.$$
For $n$ odd
$$B_p(t) = t \int^\iy_0 \psi (r)\ dr \int^1_0 (1-s)^{p} \cos\ rts\ ds.$$
For $p=0$ we have the Fourier integral formula (see \cite{Bo}, \S 9) $\:
\lim_{t\ri\iy} B_0 (t) = \pi\1 \psi(0).$
For $p \geq 1$
$$\left| t \int^\iy_1 \psi (r)\ dr \int^1_0 (1-s)^{p} \cos\ rts\ ds \right|
\leq \int^\iy_1 |\psi(r)|\ \fc{dr}{r} < \iy$$
and
$$\align
t &\int^1_0 \psi (r)\ dr \int^1_0 (1-s)^{p} \cos\ rts\ ds =
t \int^1_0 \psi (0)\ dr \int^1_0 (1-s)^{p} \cos\ rts\ ds\\ \quad\\
& + t \int^1_0 \left[\psi(r) - \psi(0)\right] dr \int^1_0 (1-s)^{p}
\cos\ rts\ ds.\endalign$$
The first integral on the right-hand side is equal to $$\psi(0)\int^1_0 (1-s)^p
s\1 \sin ts\ ds,$$ and its finiteness is well-known.
For the second integral the estimate $$\int^1_0 |\psi(r) - \psi(0)| \df{dr}{r}
< \iy$$ follows, as above. Now for $n$ even
$$\align
B_p(t) &= \int^\iy_0 \psi (r)\ r t^2 dr\ \int^1_0 (1-s)^p s J_0 (rts)\
ds\\ \quad\\
&= - \int^\iy_0 \psi' (r) \left[ rt \int^1_0 (1-s)^p s J_1 (rts)\ ds\right]
dr,\endalign$$
and integration by parts due to (7) yields the boundedness of $B_p(t)$
immediately. Therefore, the integrated terms in (10) vanish as $A\ri\iy.$ For
$\a$ integer, the element of integration coincides with that indicated in (6).

When $\a$ is fractional, apply to the first summand in the curly brackets on
the right-hand side of (10) the following formula of fractional integration by
parts (see \cite{BE2}, p. 182, or [SKM], (2.20))
$$\int^\iy_0 f_1 (t)\ R_\g (f_2\ ; t)\ dt = \int^\iy_0 W_\g (f_1\ ; t)\
f_2(t)\ dt  .$$
In our case $\g=1-\a+[\a],$
$$f_1(t)=\cases f^{([\a])}(t), & t\le A,\\ 0,& t>A,\endcases$$
and $f_2(t)={\G([\a])\over\G(\a)}{d\over dt}[t^{\a+{n\over 2}}
Q_{\a}(rt)].$ Indeed,
$$\align R_{\a-[\a]}(s^{[\a]+{n\over 2}}Q_{[\a]}(rs);t) &=
\G([\a])R_{\a-[\a]}(R_{[\a]}(s^{{n\over 2}}J_{{n\over 2}-1}(rs);\cdot);t)\\
\quad\\ =\G([\a])R_{\a}(s^{{n\over 2}}J_{{n\over 2}-1}(rs);t) &=
{\G([\a])\over\G(\a)}t^{\a+{n\over 2}}Q_{\a}(rt). \endalign$$
Therefore $f_2$ is the Riemann-Liouville derivative of order $1-\a+[\a]$
of the function $t^{[\a]+{n\over 2}}Q_{[\a]}(rt)$ and the Riemann-Liouville
integral of order $1-\a+[\a]$ of $f_2$ is exactly
$t^{[\a]+{n\over 2}}Q_{[\a]}(rt).$

Apply the usual integration by parts to the right-hand side of the formula
of the fractional integration by parts. This gives the following equality
true for all $\a:$
$$\allowdisplaybreaks\align &{1\over ([\a]-1)!}
\int^A_0 f_0^{([\a])} (t)\ t^{[\a]+\fc{n}{2}} Q_{[\a]}(rt)\,dt\\ \quad\\
&={1\over\G(\a)}\int^A_0 {d\over dt}[t^{\a+\fc{n}{2}} Q_{\a}(rt)]\,dt
{1\over \G(1-\a+[\a])}\int\limits_t^A f^{([\a])}(s)(s-t)^{[\a]-\a}\, ds\\
\quad\\ &= \fc{1}{\G(\a)\G(1-\a+[\a])}\biggl[\int^A_t (s-t)^{-\a+[\a]}
f_0^{([\a])} (s)\ ds\,\biggr] t^{\a+\fc{n}{2}} Q_\a (rt) \bigm|^A_0\\ \quad\\
&+\fc{[\a]-\a}{\G(\a)\G(1-\a+[\a])} \int^A_0 t^{\a+\fc{n}{2}} Q_\a (rt)\ dt
\int^\iy_A (s-t)^{-\a+[\a]-1} f_0^{([\a])} (s)\ ds\\ \quad \\
&-\fc{1}{\G(\a)} \int^A_0 F_\a(t)\ t^{\a+\fc{1}{2}} Q_\a (rt)\
dt.\endalign$$
It must be shown again that the two first values on the right-hand side vanish
as $A\ri\iy.$ For $t$ large enough
$$\align t^{{n\over 2}+\a} &\biggl|\int\limits_t^A (s-t)^{-\a+[\a]}
f^{([\a])}(s)\,ds \biggr| \\ \le t^{{n\over 2}+\a-[\a]} & \int\limits_t^A
(s-t)^{-\a+[\a]}|s^{[\a]}f^{([\a])}(s)|\,ds \endalign$$
and because of (2) the right-hand side equals to zero for $t=A.$
Consider now $$t^{{n\over 2}+\a}Q_{\a}(rt)\int\limits_t^A (s-t)^{-\a+[\a]}
f^{([\a])}(s)\,ds$$ as $t\to 0.$ We will use the inequality
$$Q_{\a}(rt)\le(Cr^{{n\over 2}-1})t^{{n\over 2}-1}$$ that immediately
follows from (9). As above
$$\align&\biggl| t^{{n\over 2}+\a}Q_{\a}(rt)\int\limits_1^A (s-t)^{-\a+[\a]}
f^{([\a])}(s)\,ds\biggr|\\ \quad \\ \le C&(A-t)^{1-\a+[\a]} t^{n-1+\a-[\a]}
\sup\limits_{[1,A]} |s^{[\a]} f^{([\a])}(s)| \endalign$$
and the right-hand side tends to zero as $t\to 0.$ Estimate now
$$t^{{n\over 2}+\a}Q_{\a}(rt)\int\limits_t^1 (s-t)^{-\a+[\a]}
f^{([\a])}(s)\,ds.$$
We wish to show that $$\lim\limits_{t\to 0+} t^{\a+n-1}\int\limits_t^1
(s-t)^{[\a]-\a} f_0^{([\a])}(s)\,ds=0.$$ From the definition of $F_{\a}$
we obtain $${d\over dt}\int\limits_t^\infty
(s-t)^{[\a]-\a} f_0^{([\a])}(s)\,ds=O(t^{-{n-1\over 2}}).$$ Since
$${d\over dt}\int\limits_1^\infty (s-t)^{[\a]-\a} f_0^{([\a])}(s)\,ds=
(\a-[\a])\int\limits_1^\infty (s-t)^{[\a]-\a-1} f_0^{([\a])}(s)\,ds$$
is bounded, we have
$${d\over dt}\int\limits_t^1 (s-t)^{[\a]-\a} f_0^{([\a])}(s)\,ds=
O(t^{-{n-1\over 2}})$$ and thus
$$\int\limits_t^1 (s-t)^{[\a]-\a} f_0^{([\a])}(s)\,ds=O(t^{-{n-3\over 2}}).$$
This estimate proves the statement.

Further, integrating, in $r,$ by parts $[{n-1\over 2}]$ times using
(7) with ``-'' in the left-hand side, and then once more using (7)
with ``+'' in the left-hand side, we obtain
$$\align
\int^\iy_0 &\vp_0 (r)\ r^{\fc{n}{2}} dr \int^A_0
\left\{\int^\iy_A (s-t)^{-\a+[\a]-1} f_0^{([\a])} (s)\ ds
\right\} t^{\fc{n}{2}+\a} Q_\a (rt)\ dt\\ \quad\\
&=\int^\iy_0 \psi (r)\ r^{\fc{n}{2}-\left[\fc{n+1}{2}\right]+1} dr \int^A_0
\left\{\int^\iy_A (s-t)^{-\a+[\a]-1} f_0^{([\a])} (s)\ ds \right\}\\ \quad\\
&\tm t^{\fc{n}{2}-\left[\fc{n+1}{2}\right]+\a} dt\int^1_0 (1-u)^{\a-1}
u^{\fc{n}{2}-\left[\fc{n+1}{2}\right]} J_{\fc{n}{2}+1-\left[\fc{n+1}{2}\right]}
(rtu)\ du.  \endalign$$
In view of (2) we have $\sup\limits_{s\in[A,\iy)} \left|s^{[\a]}\ f_0
^{([\a])} (s)\right| \longrightarrow O$ as $A\ri \iy.$ Besides that
$$\align
 &\int^A_0 t^{\a-[\a]-1}\ dt  \int^\iy_A (s-t)^{-\a+[\a]-1}\ ds\\ \quad\\
&= \fc{1}{\a-[\a]}  \int^1_0 t^{\a-[\a]-1} (1-t)^{-\a+[\a]}\ dt\le C.
\endalign$$
For $n$ odd and each $\a,$ since $J_{{1\over2}}(rtu)=\sqrt{{2\over \pi
rtu}}\sin rtu\,du,$ we obtain
$$\left|(rt)^{\fc{n}{2}-\left[\fc{n+1}{2}\right]+1} \int^1_0 (1-u)^{\a-1}
u^{\fc{n}{2}-\left[\fc{n+1}{2}\right]}
J_{\fc{n}{2}+1-\left[\fc{n+1}{2}\right]}\ (rtu)\ dr\right|\leq C.$$
For $n$ even Lemma 1, with $n=2,$ yields such an estimate only for
$\a\geq {1\over 2}.$
It remains to consider
$$\int^\iy_0 \psi(r)\ rdr \int^A_0\left\{\int^\iy_A
(s-t)^{-\a-1} f_0^{([\a])}\ (s)\ ds\right\} t^\a\,dt \int^1_0
(1-u)^{\a-1}\ J_1(rtu)\, du$$
for $0<\a< {1\over 2}.$
Again, applying Lemma 1 with $n=2,$ we reduce the problem to finiteness of
$$\int^\iy_1 \psi(r)\ (rt)^{\fc{1}{2} -\a} \sin rt\ dr$$
when $rt>1,$ since all the remainder terms are estimated as above after using
(7).
Let us integrate by parts.
All the integrated terms are bounded, and it remains to estimate
$$\int^\iy_1 \psi'(r)\ r(rt)^{-\fc{1}{2} -\a} \cos rt\
dr+\left(\fc{1}{2}-\a\right) \int^\iy_1 \psi(r)\ (rt)^{-\fc{1}{2} -\a} \cos rt
\ dr.$$
Both integrals are finite, and we get finally
$$(\tf,\vp)=\fc{(2\pi)^{\fc{n}{2}}(-1)^{\a^*+1}}{\G(\a)}\ \lim_{A\ri\iy}
\int^\iy_0 \vp_0(r)\ r^{\fc{n}{2}} dr \int^A_0 F_\a (t)\ t^{\a+\fc{1}{2}}
Q_\a(rt)\ dt.$$
It remains to justify the passage to the limit under the integral sign.
Due to the Lebesgue dominated theorem, it is possible if the element of
integration is dominated by an integrable function independent of $A.$
To prove this, consider two integrals: over $r\in[0,1]$ and over $r\in[1,\iy),$
respectively.
In view of (3) and (4), one can treat $F_\a$ as a function which is monotone
decreasing  and vanishing at infinity. Let $r\in[1,\iy ).$
Consider two integrals in $t\:$ over$[0,1]$ and $[1,A],$ respectively.
The first one is bounded, and the majorant is simply $C\left|\vp_0 (r)\
r^{\fc{n}{2}} \right|.$
Apply the second mean value theorem of integral calculus to the
integral over $[1,A].$ Using (7) and Lemma 1, we obtain $(\xi\leq A)\:$
$$\align
&\left|\vp_0 (r)\ r^{\fc{n}{2}}\int^A_1 F_\a (t)\ t^{\a+\fc{1}{2}} Q_\a(rt)\
dt\right|
= \left|\vp_0 (r)\ r^{\fc{n}{2}} F_\a (1)\int^\xi_1 t^{\a+\fc{1}{2}}
Q_\a(rt)\ dt\right|\\ \quad\\
=&\left| F_\a (1)\ \vp_0 (r)\ r^{\fc{n}{2}-1} \left\{ t^{\a+\fc{1}{2}}\ q_\a
(rt)\biggm|^\xi_1 +\left(\fc{n-1}{2} -\a\right)\int^\xi_1
t^{\a-\fc{1}{2}}\ q_\a (rt)\ dt\right\} \right|\\ \quad\\
\leq &\left(\fc{n-1}{2} -\a\right) \left|  F_\a (1)\ \vp_0 (r)\
r^{\fc{n-3}{2}-\a} \int^\xi_1 \fc{1}{t} \cos (rt+\mu)\ dt\right| +C \bigm|
F_\a (1)\ \vp_0 (r)\bigm| \fc{1}{r},\endalign$$
and the last value is integrable over $[1,\iy).$ At last
$$\align
&\lim_{A\ri\iy} \int^1_0 \vp_0 (r)\ r^{\fc{n}{2}} dr \int^A_0
F_\a (t)\ t^{\a+\fc{1}{2}} Q_\a(rt)\ dt\\ \quad\\
=& \lim_{A\ri\iy} \int^1_0 \left\{\left[ \vp_0 (r) -\vp_0 (0)\right] +
\vp_0 (0)\right\} r^{\fc{n}{2}}\ dt \int^A_0 F_\a (t)\ t^{\a+\fc{1}{2}}
Q_\a(rt)\ dt.\endalign$$
Since $\df{1}{r} \left[ \vp_0 (r) -\vp_0 (0)\right]$ is integrable, the part
corresponding to this function may be estimated completely like in the case
$r\in [1,\iy).$ Further,
$$\lim_{A\ri\iy} \int^1_0 r^{\fc{n}{2}} \int^A_0 F_\a (t)\ t^{\a+\fc{1}{2}}
Q_\a(rt)\ dt=\lim_{A\ri\iy}\int^A_0 F_\a (t)\ t^{\a-\fc{1}{2}}\ q_a(t)\ dt.$$
Since $F_\a$ is bounded and monotone, integration by parts and
estimates like in Lemma 1 yield the convergence of this integral in
improper sense.

In fact,we have proved the
uniform convergence of the integral (6) when $|x|\geq r_0 > 0.$

{\bf 3.2.} Let us show now that $\hf(x)\ri 0$ as $r=|x|\ri\iy.$
Using the second mean value theorem, we obtain for some $A''\leq A'\:$
$$r^{1-\fc{n}{2}} \int^{A'}_A F_\a (t)\ t^{\a+\fc{1}{2}}\ Q_\a(rt)\ dt =
F_\a (A)\ r^{-\fc{n}{2}} \left[ t^{\a+\fc{1}{2}}\ q_\a (rt)\right]^{A''}_A.$$
In view of Lemma 1, we have for every $A\geq  1$ as $A'\ri\iy$
$$\left| r^{1-\fc{n}{2}} \int^\iy_A F_\a (t)\ t^{\a+\fc{1}{2}}\ Q_\a(rt)\ dt
\right| \leq C \bigm| F_\a(A) \bigm| r^{-\fc{n+1}{2}-\a}.$$
The right-hand side tends to zero as $r\ri\iy.$ Further,
$$\align
&r^{1-\fc{n}{2}} \int^{A}_0 F_\a (t)\ t^{\a+\fc{1}{2}}\ Q_\a(rt)\ dt\\
\quad\\ & = F_\a (t)\ r^{-\fc{n}{2}} \ q_\a (rt)\bigm|^{A}_0
-r^{-\fc{n}{2}} \int^{A}_0 t^{\a+\fc{1}{2}}\ q_\a (rt)\ dF_\a (t)\\ \quad\\
&+\left(\fc{n-1}{2} -\a\right)\ r^{-\fc{n}{2}}\ \int^{A}_0 F_\a (t)\
t^{\a-\fc{1}{2}}\ q_\a (rt)\ dt = O \left(r^{-\fc{n}{2}}\right)\endalign$$
by Lemma 1 and (4).

{\bf 3.3.} Let us show the continuity of $\hf(x)$ for $|x|>0.$
Let $[r_0, r_1]$ be an interval of uniform convergence of the integral in (6),
and $|x| \in[r_0,r_1].$ Then the functions
$$\hf_k(x) = \fc{(2\pi)^{\fc{n}{2}}(-1)^{\a^*+1}}{\G(\a)}\ |x|^{1-\fc{n}{2}}
\int^k_0 F_\a (t)\ t^{\a+\fc{1}{2}}\ Q_{\a} (|x| t)\ dt$$
are continuous for each $k=1,2,\dots,\ $ and converge uniformly to $\hf(x)$ as
$k\ri\iy.$ Hence, $\hf(x)$ is continuous for these $x$ as well.

{\bf 3.4.} Let us prove now the inverse formula. Applying the Cauchy-Poisson
formula (see e.g., \cite{Bo}, Th. 56, or \cite {SW}, Ch.4, Th.3.3), we have
$$\align (2\pi)^{-n} &\int_{|u|\leq A} \left( 1-\fc{|u|^2}{A^2}\right)^
{\fc{n-1}{2}-\a}\ \hf(u)\ e^{ixu}\ du \tag11 \\ \quad\\
= \fc{(-1)^{\a^*+1}}{\G(\a)}\ r^{1-\fc{n}{2}} &\int^A_0
\left( 1-\fc{s^2}{A^2}\right)^{\fc{n-1}{2}-\a}\ sJ_{\fc{n}{2}-1}
(rs)\ ds \int^\iy_0 F_\a (t)\ t^{\a+\fc{1}{2}}\ Q_\a (st)\ dt\\ \quad\\
= \fc{(-1)^{\a^*+1}}{\G(\a)}\ r^{1-\fc{n}{2}} &\int^\iy_0 F_\a (t)\
t^{\a+\fc{1}{2}} \,dt\int^A_0 \left(1-\fc{s^2}{A^2}\right)^{\fc{n-1}{2}-\a}\
sJ_{\fc{n}{2}-1} (rs)\ Q_\a (st)\ ds.\endalign$$
This change of integrals must be justified. Let $0<\dl<A.$
The uniform convergence of the integral in $t$ for $s\geq \dl$ yields
$$\align
&\int^A_{\dl} \left( 1-\fc{s^2}{A^2}\right)^{\fc{n-1}{2}-\a}\
sJ_{\fc{n}{2}-1} (rs)\ ds \int^\iy_0 F_\a (t)\ t^{\a+\fc{1}{2}}\
Q_\a (st)\ dt\tag12\\ \quad\\ =&\int^\iy_0 F_\a (t)\ t^{\a+\fc{1}{2}}\ dt
\int^A_0 \left( 1-\fc{s^2}{A^2}\right)^{\fc{n-1}{2}-\a}\ sJ_{\fc{n}{2}-1}\
(rs)\ Q_\a (st)\ ds\\ \quad\\
-&\int^\iy_0 F_\a (t)\ t^{\a+\fc{1}{2}}\ dt
\int^\dl_0 \left( 1-\fc{s^2}{A^2}\right)^{\fc{n-1}{2}-\a}\ sJ_{\fc{n}{2}-1}
(rs)\  Q_\a (st)\ ds . \endalign$$
It suffices to show that the last integral tends to zero as $\dl\ri 0.$
Take $\ve > 0$ and let $M$ be large enough to provide $|F_\a (M)| <\ve,$
by (3). The second mean value theorem  yields after integrating by parts:
$$\align & \int^\iy_M F_\a (t)\ t^{\a+\fc{1}{2}}\ dt
\int^\dl_0 \left( 1-\fc{s^2}{A^2}\right)^{\fc{n-1}{2}-\a}\ sJ_{\fc{n}{2}-1}
(rs)\  Q_\a (st)\ ds \tag13\\ \quad\\
=F_\a (M) &\int^{M'}_M t^{\a+\fc{1}{2}}\ dt
\int^\dl_0 \left( 1-\fc{s^2}{A^2}\right)^{\fc{n-1}{2}-\a}\ sJ_{\fc{n}{2}-1}
(rs)\  Q_\a (st)\ ds\\ \quad\\    =F_\a (M) &\left[ t^{\a+\fc{1}{2}}
\int^\dl_0 \left( 1-\fc{s^2}{A^2}\right)^{\fc{n-1}{2}-\a}\ J_{\fc{n}{2}-1}
(rs)\  q_\a (st)\ ds\right]^{M'}_M\\ \quad\\
-F_\a (M) &\left(\a-\fc{n-1}{2}\right) \int^{M'}_M t^{\a-\fc{1}{2}}\ dt
\int^\dl_0 \left( 1-\fc{s^2}{A^2}\right)^{\fc{n-1}{2}-\a}\ J_{\fc{n}{2}-1}
(rs)\ q_\a (st)\ ds.\endalign$$
Estimate integrated terms in (13). The uniform boundedness,
in $t$ and $\dl,$ of the value in brackets should be shown. Since
$\z_{\a,n} \int\limits^\dl_0 \left( 1-\fc{s^2}{A^2}\right)^{\fc{n-1}{2}-\a}\
s^{-\a-\fc{1}{2}}\ J_{\fc{n}{2}-1} (rs)\ ds$ does not depend on $M$
and $M',$ such values, taking twice with opposite signs, are cancelled.
The rest, in view of (9), does not exceed, for $s\in\left[ 0,
\fc{1}{t}\right],$ the quantity
$$C t^{\a+\fc{1}{2}} \int^{\fc{1}{t}}_0 \left| J_{\fc{n}{2}-1} (rs) \right|
ds \leq C r^{\fc{n}{2}-1}\ t^{\a-\fc{n-1}{2}},$$
and for $t\dl > 1$ and $s\in\left[{1\over t}, \dl\right],$ in view of Lemma 1
and (9), is
$$\G(\a) t^{{1\over 2}} \int^\dl_{\fc{1}{t}}  s^{-\a}
\left( 1-\fc{s^2}{A^2}\right)^{\fc{n-1}{2}-\a}\ J_{\fc{n}{2}-1}\ (rs)
J_{\fc{n}{2}+\a}
\ (st)\ ds + O\left(\fc{1}{t}\int^\dl_{\fc{1}{t}} s^{\fc{n-5}{2}-\a}
ds\right).$$
The remainder term does not exceed $C\left\{ t^{\a-\fc{n-1}{2}} + (t\dl)\1
\dl^{\fc{n-1}{2}-\a}\right\},$ which is bounded.
Apply (9) to $J_{\fc{n}{2}-1}$ and (8) to $J_{\fc{n}{2}+\a}$ in the main
term. We obtain for $\a < {n-1\over 2}$
$$\int^\dl_{\fc{1}{t}} \left( 1-\fc{s^2}{A^2}\right)^{\fc{n-1}{2}-\a}
s^{-\a+\fc{n}{2}-\fc{3}{2}} ds \leq C\left( \dl^{\fc{n-1}{2}-\a} +
t^{\a-\fc{n-1}{2}}\right) \leq C.$$
For $\a= {n-1\over 2},$ integrate by parts as follows:
$$\align
&t^{\fc{1}{2}} \int^\dl_{\fc{1}{t}} s^{-\fc{n-1}{2}}\ J_{\fc{n}{2}-1}\
(rs)\ J_{n-\fc{1}{2}}\ (st)\ ds\\ \quad\\
&= - t^{-\fc{1}{2}} s^{-\fc{n-1}{2}}\ J_{\fc{n}{2}-1}\ (rs)\
J_{n+\fc{1}{2}}\ (st) \biggm|^\dl_{\fc{1}{t}}
+ rt^{-\fc{1}{2}}\int^\dl_{\fc{1}{t}} s^{-\fc{n-1}{2}}
 J_{\fc{n}{2}-2}\ (rs)\ J_{n+\fc{1}{2}}\ (st)\ ds,\endalign$$
and now the proof is continued as for the remainder term.

Let us estimate the last integral in (13). This makes sense only for $\a <
{n-1\over 2}.$ Use again Lemma 1. We have to estimate
$$\int^{M'}_M t^{-\fc{1}{2}} dt \int^\dl_{\fc{1}{t}}
\left( 1-\fc{s^2}{A^2}\right)^{\fc{n-1}{2}-\a} s^{-\a}  J_{\fc{n}{2}-1}\ (rs)\
J_{\fc{n}{2}+\a}\ (ts)\ ds,$$
since the rest is estimated analogously. Let us change the order of
integration. Without loss of generality, one can take $\delta <
{1\over  M}.$  The following should be estimated:
$$\align
& \int^\dl_{\fc{1}{M'}} \left( 1-\fc{s^2}{A^2}\right)^{\fc{n-1}{2}-\a} s^{-\a}
J_{\fc{n}{2}-1}\ (rs)\ ds \int_{\fc{1}{s}}^{M'} t^{-\fc{1}{2}}
J_{\fc{n}{2}+\a}\ (st)\ dt\\ \quad\\       +&\int_\dl^{\fc{1}{M}}
\left( 1-\fc{s^2}{A^2}\right)^{\fc{n-1}{2}-\a} s^{-\a}  J_{\fc{n}{2}-1}\
(rs)\ ds \int_M^{\fc{1}{s}} t^{-\fc{1}{2}} J_{\fc{n}{2}+\a}\ (st)\
dt.\endalign$$
Apply (8) to $J_{\fc{n}{2}+\a}.$
For the remainder term, after applying (9) to $J_{\fc{n}{2}-1},$ we obtain
$$\int^\dl_{\fc{1}{M'}} s^{\fc{n-5}{2}-\a} ds \int_{\fc{1}{s}}^{M'}
\fc{dt}{t^2} + \int_\dl^{\fc{1}{M}} s^{\fc{n-5}{2}-\a} ds \int_M^{\fc{1}{s}}
\fc{dt}{t^2}$$
which is obviously bounded. For the main term, the integral in $t$ is of the
form $$\int t^{-1} \cos (st+\mu)\ dt.$$ After integrating by parts the
estimates coincide with those for the remainder term.
Hence, the value (13) is small. Choosing $\dl$ so small that
$$\left| \int^M_0 F_\a (t)\ t^{\a+\fc{1}{2}} dt \int^\dl_0
\left( 1-\fc{s^2}{A^2}\right)^{\fc{n-1}{2}-\a} s J_{\fc{n}{2}-1} (rs)\ Q_\a
(st)\ ds \right| < \ve,$$
we obtain that the last integral in (12) tends to zero as $\dl\ri 0.$
Returning to (11) we have
$$\align
&\int^A_0 \left( 1-\fc{s^2}{A^2}\right)^{\fc{n-1}{2}-\a} s
J_{\fc{n}{2}-1} (rs)\ Q_\a (st)\ ds = \int^A_0 s J_{\fc{n}{2}-1}
(rs) Q_\a (st)\ ds\\ \quad\\
-&\int^A_0 s J_{\fc{n}{2}-1} (rs)\ Q_\a (st)\ ds\ 2\left({\fc{n-1}{2}-\a}
\right) \int_0^{\fc{s}{A}} u (1-u^2)^{\fc{n-3}{2}-\a} du.\endalign$$
Let us begin with the second integral. We have
$$\align
&\int^\iy_0 F_\a (t)\ t^{\a+\fc{1}{2}} dt \int^A_0 s J_{\fc{n}{2}-1} (rs)\
Q_\a (st)\ ds \int_0^{\fc{s}{A}} u (1-u^2)^{\fc{n-3}{2}-\a} du\\ \quad\\
=&\int^1_0 u (1-u^2)^{\fc{n-3}{2}-\a} du \int^\iy_0 F_\a (t)\
t^{\a+\fc{1}{2}} dt \int^A_{Au} s J_{\fc{n}{2}-1} (rs)\ Q_\a (st)\ ds\\
\quad\\ =&\int^1_0 u (1-u^2)^{\fc{n-3}{2}-\a} du \left\{ F_\a (t)\
t^{\a+\fc{1}{2}} \int^A_{Au} J_{\fc{n}{2}-1} (rs)\ q_\a (st)\
ds\biggm|^\iy_0\right.\\ \quad\\
-&\int^\iy_0 d F_\a (t)\ t^{\a+\fc{1}{2}} \int^A_{Au} J_{\fc{n}{2}-1}
(rs)\ q_\a (st)\ ds\\ \quad\\
+&\left.\left({\fc{n-1}{2}-\a}\right) \int^\iy_0  F_\a (t)\
t^{\a-\fc{1}{2}} dt \int^A_{Au} J_{\fc{n}{2}-1} (rs)\ q_\a (st)\
ds\right\}.\endalign$$
The integrated terms for $t=0$ and $t=\infty$ vanish. Indeed, it is
obvious for $t=0$ and for $t=\infty$ follows from Lemma 1 and (3).
We have to show that the right-hand side tends to zero as $A\ri\iy.$
Note firstly, that the estimate $|q_\a (st)| \leq C (st)^{-\a-\fc{1}{2}}$ and
estimates (8) and (9) for $J_{\fc{n}{2}-1}$ yield
$$\left| t^{\a+\fc{1}{2}} \int^A_{Au} J_{\fc{n}{2}-1} (rs)\ q_\a
(st)\ ds \right| \leq C \int^A_{Au} \fc{ds}{s^{1+\ve}}$$
with some $\ve\in (0, 1).$
This estimate combined with (4) gives proper estimates for
the second summand in the curly brackets.
When $t\in [0,1]$ the estimates for the third one are similar.
Use then Lemma 1 when $t\in[1,\iy).$ The remainder term is estimated as above.
The main term, by the second mean value theorem, is equal to
$$\align
&\int^\iy_1 F_\a (t)\ t^{-\fc{1}{2}} dt\int^A_{Au} s^{-\a} J_{\fc{n}{2}-1}
(rs) J_{\fc{n}{2}+\a} (st)\ ds\\ \quad\\
= F_\a (1)& \int^\xi_1 t^{-\fc{1}{2}} dt\int^A_{Au} s^{-\a} J_{\fc{n}{2}-1}
(rs) J_{\fc{n}{2}+\a} (st)\ ds.\endalign$$
Apply (8) to $J_{\fc{n}{2}+\a}.$
The remainder term does not need new techniques. It remains to estimate
$$\int^\xi_1 \fc{dt}{t} \int^A_{Au} s^{-\a-\fc{1}{2}} J_{\fc{n}{2}-1} (rs)
\cos (ts+\mu)\ ds.$$
Integration by parts in $t$ and estimates of the integral in $s$, like above,
finish the proof of a tendency of the limit to zero. It remains to consider
$$\align
&\int^\iy_0 F_\a (t)\ t^{\a+\fc{1}{2}} dt \int^A_0 s J_{\fc{n}{2}-1}
(rs)\ Q_\a (st)\ ds\\ \quad\\
=&\int^\iy_0 F_\a (t)\ t^{\a-\fc{n-1}{2}} \fc{d}{dt} \left[ t^{\fc{n}{2}}
 \int^1_0 (1-u)^{\a-1} u^{\fc{n}{2}-1} du \int^A_0 J_{\fc{n}{2}-1} (rs)
J_{\fc{n}{2}} (uts)\ ds\right] dt.\endalign$$
Let us substitute integration over $[0,A]$ by integration over the difference
of two sets: $[0,\iy)$ and $[A,\iy).$  Let us use the formula
$$\G (\nu-\mu) \int^\iy_0 J_\mu (at) J_\nu (bt) t^{\mu-\nu+1} dt=\cases
2^{\mu-\nu+1} a^\mu b^{-\nu} (b^2-a^2)^{\nu-\mu+1}, &\text{for} \ \ b>
a,\\ \quad \\   0, &\text{for} \ \ b<a\endcases$$
which is true for $\nu>\mu>-1$ (see \cite{BE2}, p. 148). We obtain
$$\align
&\int^\iy_0 F_\a (t)\ t^{\a-\fc{n-1}{2}} \fc{d}{dt} \left[ t^{\fc{n}{2}}
 \int^1_0 (1-u)^{\a-1} u^{\fc{n}{2}-1} du \int^\iy_0 J_{\fc{n}{2}-1} (rs)
 J_{\fc{n}{2}} (uts)\ ds\right] dt\\ \quad\\
=&\int^\iy_r F_\a (t)\ t^{\a-\fc{n-1}{2}} \fc{d}{dt} \left[ t^{\fc{n}{2}}
\int^1_{\fc{r}{t}} (1-u)^{\a-1} u^{\fc{n}{2}-1} r^{\fc{n}{2}-1}
(tu)^{-\fc{n}{2}} du\right] dt\\ \quad\\
=r^{\fc{n}{2}-1}& \int^\iy_r f_0^{(\a)} (t)\ (t-r)^{\a-1} dt.\endalign$$
Integrating by parts and using (1)--(3), we get
$$\fc{(-1)^{\a^*+1}}{\Gamma(a)}\ \int^\iy_r f_0^{(\a)} (t)\ (t-r)^{\a-1} dt =
f_0(r).$$ For $\a$ fractional we used also the permutability of the
fractional integral and fractional derivative. If we show that
$$\int^\iy_0 F_\a (t)\ t^{\a-\fc{n-1}{2}} \fc{d}{dt}
\left[ t^{\fc{n}{2}} \int^1_0 (1-u)^{\a-1} u^{\fc{n}{2}-1} du \int^\iy_A
J_{\fc{n}{2}-1} (rs) J_{\fc{n}{2}} (uts)\ ds\right] dt$$
tends to zero as $A\ri\iy$ the inverse formula will be proved.
Integration by parts in the outer integral yields
$$\align &F_\a (t)\ t^{\a+\fc{1}{2}}
 \int^1_0 (1-u)^{\a-1} u^{\fc{n}{2}-1} du \int^\iy_A J_{\fc{n}{2}-1} (rs)
 J_{\fc{n}{2}} (ust)\ ds\biggm|^\iy_0\\ \quad\\
&-\int^\iy_0 \left[ t^{\a+\fc{1}{2}}
 \int^1_0 (1-u)^{\a-1} u^{\fc{n}{2}-1} du \int^\iy_A J_{\fc{n}{2}-1} (rs)
 J_{\fc{n}{2}} (ust)\ ds \right] d F_\a (t)\\ \quad\\
&+ \left(\fc{n-1}{2} - \a\right) \int^\iy_0 F_\a (t)\
t^{\a-\fc{1}{2}} dt  \int^1_0 (1-u)^{\a-1} u^{\fc{n}{2}-1}
du \int^\iy_A J_{\fc{n}{2}-1} (rs)  J_{\fc{n}{2}} (ust)\ ds.\endalign$$
To estimate the first two summands, in view of (4), it suffices to show that
$$\lim_{A\ri\iy} \sup_{t} \left| t^{\a+ \fc{1}{2}} \int^1_0 (1-u)^{\a-1}
u^{\fc{n}{2}-1} du \int^\iy_A J_{\fc{n}{2}-1} (rs)  J_{\fc{n}{2}} (ust)\ ds
\right| = 0.$$
Let us show that the change of the order of integration is legal.
By the Lebesgue dominated theorem, it suffices to find
an integrable majorant, independent of $A$ and $A',$ for the function
$$(1-u)^{\a-1} u^{\fc{n}{2}-1} \int^{A'}_A J_{\fc{n}{2}-1} (rs)
 J_{\fc{n}{2}} (uts)\ ds. $$
Integrating by parts in the inner integral, we obtain
$$\align & \qquad\qquad\qquad
 \int^{A'}_A J_{\fc{n}{2}-1} (rs)\  J_{\fc{n}{2}} (uts)\ ds \\ \quad\\
& = \fc{1}{r}\    J_{\fc{n}{2}} (rs)\ J_{\fc{n}{2}} (uts) \biggm|^{A'}_A
+ \fc{tu}{r} \int^{A'}_A J_{\fc{n}{2}} (rs)\ J_{\fc{n}{2}+1} (uts)\
ds.\endalign$$
The integrated terms are bounded, and $(1-u)^{\a-1} u^{\fc{n}{2}-1}$ is
integrable over $[0, 1].$ Apply (8) to the Bessel functions in the last
integral. The remainder terms are bounded, and $(1-u)^{\a-1} u^{\fc{n-3}{2}}$
is integrable. We have for the main terms
$$\align
&\int^{A'}_A \fc{\cos (rs-\fc{\pi n}{2}-\fc{\pi}{4}) \cos
(uts-\fc{\pi n}{2}     -\fc{3\pi}{4})}{s}\ ds\\ \quad\\
=\fc{1}{2}\ &\int^{A'}_A \fc{\cos (rs+uts -\pi n -
\pi)}{s}\ ds - \fc{1}{2}\ \int^{A'}_A \fc{\sin s (r-tu)}{s}\ ds,\endalign$$
and the boundedness of these integrals is easily obtained by integration
by parts and inequality $(r+ut)\1\leq r\1.$
Now Lemma 1 and the rough estimate $J_\nu (t) = O(t^{-\fc{1}{2}})$ yield
$$\align
&\lim_{A\ri\iy} \sup_{t} \left| t^{\a+ \fc{1}{2}} \int^1_0 (1-u)^{\a-1}
u^{\fc{n}{2}-1} du \int^\iy_A J_{\fc{n}{2}-1} (rs)\  J_{\fc{n}{2}} (uts)\ ds
\right| \\ \quad\\
=&\lim_{A\ri\iy} \sup_{t} \left| t^{\a+ \fc{1}{2}} \int^\iy_A
 J_{\fc{n}{2}-1} (rs)\ q_\a(ts)\ ds \right| \leq C \lim_{A\ri\iy}
r^{-\fc{1}{2}} \int^\iy_A \fc{ds}{s^{1+\a}}= 0 .\endalign$$
It remains to estimate
$$\align
&\int^\iy_0 F_\a (t)\ t^{\a-\fc{1}{2}} dt
 \int^1_0 (1-u)^{\a-1} u^{\fc{n}{2}-1} du \int^\iy_A J_{\fc{n}{2}-1}
(rs)\ J_{\fc{n}{2}} (uts)\ ds\\ \quad\\
= &\int^\iy_0 F_\a (t)\ t^{\a-\fc{1}{2}} dt \int^\iy_A J_{\fc{n}{2}-1}
(rs)\ q_\a (ts)\ ds.\endalign$$
For $t\in[{r\over 2},\iy]$ such estimates are already made above.
Let $t\in[0,{r\over 2}].$ Consider
$$ \int\limits_0^{{r\over 2}}F_{\a}(t)t^{\a-{1\over 2}}\,dt
\int\limits_0^1 (1-u)^{\a-1}u^{{n\over 2}-1}\,du\int\limits_A^\infty
J_{{n\over 2}-1}(rs)J_{{n\over 2}}(ust)\,ds .$$
Apply (8) to the first Bessel function on the right-hand side.
Estimates for the remainder terms are obvious, so we have to estimate
$$\int\limits_0^{{r\over 2}}F_{\a}(t)t^{\a-{1\over 2}}\,dt
\int\limits_0^1 (1-u)^{\a-1}u^{{n\over 2}-1}\,du\int\limits_A^\infty
{\cos(rs-h)\over \sq{s}}J_{{n\over 2}}(ust)\,ds$$ where $h$ is a number.
Integration by parts in the integral $$\int\limits_A^\infty
{\cos(rs-h)\over \sq{s}}J_{{n\over 2}}(ust)\,ds=\int\limits_A^\infty
{\cos(rs-h)\over \sq{s}s^{{n\over 2}}}s^{{n\over 2}}J_{{n\over 2}}(ust)\,ds$$
and simple calculations using (7) and (8) yield that the problem
is to show that the integral
$$ut\int\limits_A^\infty J_{{n\over 2}-1}(ust){\sin(rs-h)\over\sq{s}}\,ds$$
tends to zero as $A\to \infty.$ We are going to apply (8) to the
Bessel function. For this we had to get the factor $ut$ before the integral.
It provides convergence of other integrals, in $t$ and $u.$ Again for the
remainder terms in (8) estimates are similar to those afore-mentioned.
And for the main term we get the inner integral in the form
$$\sq{ut}\int\limits_A^\infty{\cos(ust-l)\sin(rs-h)\over s}\,ds.$$
Since $ut\le{r\over 2}$ for $t\in[0,{r\over 2}]$ and $u\in[0,1],$
the latter integral is obviously small as $A\to\infty.$
The proof is complete.\ \ \ \qed
\bk
\subheading{4. FOURIER TRANSFORM: \newline
FROM MANY DIMENSIONS TO ONE DIMENSION}
\bk
{\bf 4.1.} It happens sometimes that certain conditions provide that the
difference is bounded between an integral of the absolute value of the
multiple Fourier transform and an integral of the absolute value of the
one-dimensional Fourier transform of some other function.
Such a result was obtained by Podkorytov \cite{P}, provides that the function
under consideration is radial and boundedly supported, and its Fourier
transform is integrable. In the following theorem this process is realized
under essentially different conditions. This will allow us to generalize some
important one-dimensional results to the multi-dimensional case.
Besides that, we get an easier way to estimate the growth of the Fourier
transform when it is non-integrable.
\sk
\proclaim{Theorem 2}\ Let a function $f$ be radial and satisfying conditions
{\rm (1) - (4)} with $\a={n-1\over 2}.$ Suppose further that
$$\int^1_0 \fc{| F(t) |}{t}\ dt < \iy.\tag14$$      Then for $n\geq 2$
$$\align \int_{1\leq|x|\leq N} |\hf(x)| dx &= \fc{2^{\fc{n+3}{2}}
\pi^n}{\G\left(\fc{n}{2}\right)} \int^N_1\left| \int^\iy_0 F(t)\ \sin
(rt-\fc{\pi n}{2})\ dt\right| dr\\ \quad\\
&+ \th \left(V_F +  \int^1_0 \fc{|F(t)|}{t}\ dt\right).\endalign$$
where $|\th | \leq C.$\ep
\sk
\demo{Proof}\ After using (6) and passage to spherical coordinates,
multidimensional integration is reduced to the one-dimensional,
and one has to estimate
$$\int^N_1 r^{\fc{n}{2}} \left| \int^\iy_0 F(t)\  t^{\fc{n}{2}} Q (rt)\
dt\ \right| dr.$$
When $rt < 1,$ we get, as when proving Lemma 2
$$(rt)^{\fc{n}{2}} Q(rt) = \G\left(\fc{n-1}{2}\right) (rt)^{\fc{1}{2}}
J_{n-\fc{3}{2}} (rt) + O (1).$$
Now (9) yields
$$\align
&\int_{1\leq|x|\leq N} |\hf(x)| dx =  \int^N_1 \left| \fc{2^{\fc{n+3}{2}}
\pi^n}{\G\left(\fc{n}{2}\right)} \int^{\fc{1}{r}}_0 F(t)\ \sin
\left(rt-\fc{\pi n}{2}\right)\ dt\right.\\ \quad\\
&+ \left.r^{\fc{n}{2}} \int_{\fc{1}{r}}^\iy F(t)\ t^{\fc{n}{2}} Q (rt)\ dt
\right| dr + \th \int^N_1 \int^{\fc{1}{r}}_0 |F(t)|\ dt dr.\endalign$$
Observe that
$$\int^N_1 \int^{\fc{1}{r}}_0 |F(t)|\ dt dr \leq    \int^\iy_1
\int^{\fc{1}{r}}_0 |F(t)|\ dt dr \leq \int_0^1\fc{|F(t)|}{t}\ dt .$$
Take now $rt \geq 1.$
Applying Lemma 2 with $p=2, \mu=\fc{n-3}{2},\ \be= \fc{n}{2} - 1,$
we get
$$\align
(rt)^{\fc{n}{2}} Q(rt) &= \G\left(\fc{n-1}{2}\right) (rt)^{\fc{1}{2}}
J_{n-\fc{3}{2}} (rt)\\ \quad\\
&+ \fc{n-3}{2}\ \G\left(\fc{n+1}{2}\right) (rt)^{-\fc{1}{2}}
J_{n-\fc{1}{2}} (rt) + \th (rt)\2.\endalign$$
The remainder is estimated easily by
$$\align
&\int^N_1 r\2 dr \int_{\fc{1}{r}}^\iy |F(t)|\ t\2 dt \leq \int^1_0 |F(t)|\ t\2
dt \int_{\fc{1}{t}}^\iy r\2 dr\\ \quad\\
&+ \int_1^\iy |F(t)|\ t\2 dt \int_1^\iy r\2 dr\leq \int^1_0 |F(t)| \fc{dt}{t}
+ V_F.\endalign$$
Take away the main term of $\G\left(\fc{n-1}{2}\right) (rt)^{\fc{1}{2}}
J_{n-\fc{3}{2}} (rt)$ (see (8)) from the asymptotic expansion for
$(rt)^{{n\over 2}}Q(rt)$ and estimate the integral of the function $F$
multiplied by
$$\align & \G\left(\fc{n-1}{2}\right) (rt)^{\fc{1}{2}}\left[J_{n-\fc{3}{2}}
(rt) +\sqrt{\fc{2}{\pi rt}} \sin \left(rt - \fc{\pi n}{2}\right)\right] \\
\quad\\  &+\fc{n-3}{2} \G \left(\fc{n+1}{2}\right) (rt)^{-\fc{1}{2}}
J_{n-\fc{1}{2}} (rt). \endalign$$
Integrating by parts, we obtain
$$\align
&\int^N_1 r^{\fc{1}{2}} \left| F(t) \int_t^\iy s^{\fc{1}{2}}
\left[ J_{n-\fc{3}{2}} (rs) +\sqrt{\fc{2}{\pi rs}} \sin \left(rs - \fc{\pi
n}{2}\right)\right] ds \right. \biggm|^\iy_{\fc{1}{r}}\\ \quad\\
&\left.-\int^\iy_{\fc{1}{r}} dF(t) \int^\iy_t s^{\fc{1}{2}}
\left[ J_{n-\fc{3}{2}} (rs) +\sqrt{\fc{2}{\pi rs}} \sin \left(rs - \fc{\pi
n}{2}\right)\right] ds \right|\ dr\\ \quad\\
&+ \int^N_1 r^{-\fc{3}{2}} \left| F(t)\ t^{-\fc{1}{2}} J_{n-\fc{3}{2}} (rt)
\bigm|^\iy_{\fc{1}{r}} + \int^\iy_{\fc{1}{r}} t^{-\fc{1}{2}}
J_{n-\fc{3}{2}} (rt)\ dF (t)\right.\\ \quad\\
&\left.-\fc{1}{2}\ \int^\iy_{\fc{1}{r}} F(t)\ t^{-\fc{3}{2}} J_{n-\fc{3}{2}}
(rt)\, dt \right| \ dr.\endalign$$
Taking into account $\left|J_{n-\fc{3}{2}} (rt) \right|
\leq C(rt)^{-\fc{1}{2}},$ the second integral (in $r$) is estimated by
$$\align
&\int^N_1 \left| F\left(\fc{1}{r}\right)\right|\ \fc{dr}{r} + \int^N_1 r\2
dr \int^\iy_{\fc{1}{r}} \fc{1}{t} \left| dF(t)\right|\\ \quad\\
+ &\int^N_1 r\2 dr \int^\iy_{\fc{1}{r}} |F(t)| t\2 dt \leq 2 \int^1_0
|F(t)|\ \fc{dt}{t}\\ \quad\\
+ V_F + &\int^1_0 \fc{1}{t} | dF(t)| \int^\iy_{\fc{1}{t}} r\2 dr
+ \int^\iy_1 \fc{1}{t} | dF(t)| \int^\iy_1 r\2 dr\\ \quad\\
\leq 2&\left( V_F + \int^1_0 |F(t)|\ \fc{dt}{t}\right).\endalign$$
The integrated term under the sign of absolute value in the first integral in
$r$ obviously vanishes at infinity. For $t=\df{1}{r}$ use (8) and obtain
$$\align
&\left| \int^\iy_{\fc{1}{r}} s^{\fc{1}{2}} \left[
J_{n-\fc{1}{2}} (rs) -\sqrt{\fc{2}{\pi rs}} \sin \left(rs - \fc{\pi n}{2}
\right)\right] ds \right|\\ \quad\\
\leq C r^{-\fc{3}{2}} &\left| \int^\iy_{\fc{1}{r}} \fc{1}{s} \sin
\left(rs - \fc{\pi n}{2}\right)\ ds\right| + C r^{-\fc{5}{2}} \left|
\int^\iy_{\fc{1}{r}} s\2 ds\right| = O( r^{-\fc{3}{2}}),\endalign$$
and
$$\int^N_1\left| F\left(\fc{1}{r}\right)\right| \fc{dr}{r} \leq
\int^1_0 |F(t) |\ \fc{dt}{t}.$$
Apply again (8) to the integral remained and get
$$\align
&\int^N_1 \fc{dr}{r} \left| \int^\iy_{\fc{1}{r}} dF(t) \int^\iy_t
\fc{1}{s} \sin \left(rs - \fc{\pi n}{2}\right)\ ds \right|\\ \quad\\
+& \int^N_1 r\2 dr \int^\iy_{\fc{1}{r}} |dF(t) | \int^\iy_t s\2 ds
\leq 3 \int^\iy_1 r\2 dr \int^\iy_{\fc{1}{r}} \fc{1}{t} |dF(t)|.
\endalign$$
But the right-hand side was already estimated above.
The proof of the theorem is complete. \ \ \ \qed\edm
\sk \sk
\remark{Remark 3}\  Condition (14) is essential. Let $n=3.$
Consider the function
$$f(x)=\sin\biggl(\ln\ln{e\over|x|}\biggr)$$
for $|x|\in[0,1],$ and $0$ otherwise. We have
$$f_0 '(t)={1\over t\ln{e\over t}\ln\ln{e\over t}}\cos
\biggl(\ln\ln{e\over t}\biggr)$$
and $F(t)=tf_0 '(t).$ This function obviously satisfies conditions (1)-(4).
It is easy to see that
$$\int^1_0 {|F(t)|\over t}\,dt=\int^1_0|f'(t)|\,dt=\iy.$$ Consider
$$\align  \int\limits_{1\le|x|\le N} |\hat f(x)|\,dx &=C\int^N_1
r^{{3\over 2}}\,dr\biggl|\int_0^1 F(t)t^{{3\over 2}}\,dt
\int_0^1 s^{{3\over 2}} J_{{1\over 2}} (rts)\,ds\biggr| \\ \quad\\
&=C\sqrt{{2\over\pi}}\int_1^N rdr\biggl|\int_0^1 F(t) tdt\int_0^1
s\sin rts\,ds\biggr| \\ \quad\\
&=C\sqrt{{2\over\pi}}\int_1^N dr\biggl|\int_0^1 F(t)
\cos rt\,dt -\int_0^1 {1\over rt}F(t)\sin rt\,dt\biggr|.\endalign$$
It suffices now to prove that
$$\lim\limits_{N\to\iy}\int_1^N{dr\over r}\biggl|
\int_0^1 {F(t)\over t}\sin rt\,dt\biggr|=\iy.$$
We have after integrating by parts
$$\int_1^N{dr\over r}\biggl|\int_0^1 {F(t)\over t}\sin rt\,dt\biggr|=
\int_1^N dr \biggl|\int_0^1 f_0(t)\cos rt\,dt\biggr|,$$
and one has to prove that the one-dimensional Fourier transform of the
function $f_0$ is nonintegrable. Indeed, if it were integrable, the
following condition would valid necessarily (see e.g., [K], Ch.2, \S 10):
the integral $\int_0^1 {f_0 (t)\over t}\,dt$ converges. But it is easy to
see that this integral diverges for our function. Thus, we have built the
counterexample which shows that condition (14) is sharp.

Another, somewhat more complicate example suggested by E. Belinskii was
given in [L2]:
$$f(x)={1\over\ln\ln{e\over|x|}}\sin\biggl(\ln\ln{e\over|x|}\biggr).$$
It was constructed in order $f$ to be continuous at the origin.
\endremark
\sk
{\bf 4.2.} Let us obtain a generalization, to the multiple case, of the
Zygmund-Bochkarev criterion for the absolute convergence of Fourier
series of the function of bounded variation (see \cite{B1},
\cite{B2}, Ch.2, Th.3.1).
Let us note that an integral analog of the Stechkin criterion (see e.g.
\cite{Ba}) is proved by Trigub (\cite{T4}, Th. 2). We use the standard notation
$\omega$ for the modulus of continuity.
\sk
\proclaim{Corollary 1}\ Let a radial function be boundedly supported,
satisfy conditions {\rm (1)} and {\rm (4)}, and $F(0) = 0.$ Then the condition
$$\sum^\iy_{k=1} \fc{1}{k} \sqrt{\om \left( F\ ; \fc{1}{k}\right)} <
\iy\tag15$$
is sufficient and, on the whole class, necessary for $\hf\in L^1 (\BR^n).$\ep
\sk
\demo{Proof}\ It is natural to consider those functions which are satisfying
$\sum\limits^\iy_{k=1} \fc{1}{k}\ \om \left( F ;\fc{1}{k}\right) <\iy.$
Since $|F(t)| \leq \om \left( F\ ;\fc{1}{k}\right)$, the condition
$\sum\limits^\iy_{k=1} \fc{1}{k} \left| F(\fc{1}{k})\right| <\iy$
provides that (14) holds.
Now we are able to apply Theorem 2 and to reduce the problem to the
one-dimensional one (see \cite{B1} or \cite{B2}, Ch.2, \S 3, Th.3.1).
True, the absolute convergence of the Fourier series was investigated there,
but it is closely connected with the integrability of the Fourier transforms
due to the following theorem of Trigub (see \cite{T0}, \cite{T1};
an application to summability is given in [T2], Corollary 2):
\sk
\proclaim{Theorem A1}\ Let $f(t)$ be a boundedly supported
function of one variable, and $f_1(t) = tf(t).$
Then $\hf\in L^1 (\BR^1)$ if and only if the functions $f$ and
$f_1$ after periodic continuation have absolutely convergent
Fourier series. \endproclaim \flushpar
It remains to note that $\om \left( F\ ; \fc{1}{k}\right) \leq C\max
\{ \om (F,t)\ ; t\},$ and the corollary is proved.  \qed \enddemo
\sk

{\bf 4.3.} Let us give two examples.
\example{Example 1} \ Consider the function $f(x) =(1-|x|^\a)^\be_+.$
It was proved earlier (see \cite{Lf}, \cite{T3}) that for $\a >0$ and $\be >
\fc{n-1}{2}$ one has $\hf\in L^1 (\BR^n).$
Let us establish this fact by means of Corollary 1.
Condition (1) is evidently satisfied.
The same may be said about conditions (4) and (15) for $n$ odd.
For $n$ even to verify (4), it suffices to show that for $\psi(t)
= t^\g (a-t^\a)^{\fc{1}{2} + \ve}_+,$ with $\ve \geq 0,\ \g \geq 0,$
the function $t^{\fc{1}{2}} \psi^{(\fc{1}{2})} (t)$ is of bounded
variation. We have
$$\align
\psi^{(\fc{1}{2})} (t) &= \fc{d}{dt} \int^1_t (s-t)^{-\fc{1}{2}}\ s^\g (1-s^\a)
^{\fc{1}{2}+\ve}\ ds\\ \quad\\
&= \int^1_t (s-t)^{-\fc{1}{2}}\ \fc{d}{ds} \left\{ s^\g (1-s^\a)
^{\fc{1}{2}+\ve}\right\}\ ds.\endalign$$
This means that the boundedness of variation of the function $$t^{\fc{1}{2}}
\int^1_t (s-t)^{-\fc{1}{2}}\ s^{\z-1} (1-s^\a)^{\ve-\fc{1}{2}} ds,$$
with $\ve\geq 0,\ \z> 0,$ should be established. Further,
$$s^{\z-1} (1-s^\a)^{\ve-\fc{1}{2}} =(1-s^\a)^{\ve-\fc{1}{2}} (s^{\z-1} - 1) +
(1-s^\a)^{\ve-\fc{1}{2}}.$$
Denoting $C_\a = \lim\limits_{s\ri 1}
\left(\fc{1-s^\a}{1-s}\right)^{\ve-\fc{1}{2}},$  we have  \mk
$$\align
(1-s^\a)^{\ve-\fc{1}{2}} &= (1-s^\a)^{\ve-\fc{1}{2}} -
C_\a (1-s)^{\ve-\fc{1}{2}} + C_\a(1-s)^{\ve-\fc{1}{2}}\\ \quad\\
&= (1-s^\a)^{\ve+\fc{1}{2}} \fc{\left(\fc{1-s^\a}{1-s}\right)^{\ve-\fc{1}{2}}
-C_\a}{1-s} + C_\a (1-s)^{\ve-\fc{1}{2}}.\endalign$$ \mk\flushpar
Thus, the boundedness of variation of the function
$$t^{\fc{1}{2}} \int^1_t (s-t)^{-\fc{1}{2}}\ (1-s)^{\ve-\fc{1}{2}}\ ds$$
should be established.
But, after a simple change of variables, it is equal to the function
$$t^{\fc{1}{2}} (1-t)^\ve \int^1_0 s^{-\fc{1}{2}}\
(1-s)^{\ve-\fc{1}{2}}\ ds,$$ and (4) is now obvious.
Moreover, this makes (14) obvious too, and this completes the proof of the
example.\endexample
\sk
\example{Example 2}\ Consider $f(x) = \df{1-(1-|x|^\a)^\be_+}{|x|^r}$
with $\a > r$ and $\be > \df{n-1}{2}.$
Let us show again that $\hf\in L^1 (\BR^n).$
On $[0,1]$ the argument from Example 1 is applicable.
Hence for $t\in [0,1]$ we have
$$F(t) = C_1 t^{\a-r} (1-t)^{\be -\fc{n-1}{2}} + g(t),$$
where $g$ is a continuously differentiable function, $g(0)=0.$
It is clear that (14) is satisfied.
Applying, if needed, the formula of fractional
derivation (see \cite{BE2}, p. 201), we obtain
$t^{\fc{n-1}{2}} \left(t^{-r}\right)^{(\fc{n-1}{2})} = C_2 t^{-r}.$
Use now Theorem 2 and integrate by parts in the one-dimensional integral.
After that we have to estimate the following value:
$$\int^N_1 \fc{1}{s} \left| \int^\iy_0 F'(t) \cos
\left(st - \fc{\pi n}{2}\right)\ dt\right|\ ds.$$
On $[0,1]$ we have $F'(t) \in \Lip \ve$ in $L^1$ metric, for some $\ve > 0.$
But the integral $\int^N_1 \df{ds}{s^{1+\ve}}$ converges.
For $t\in [1,\iy],$ we obtain by integration by parts:  \mk
$$\align
&\int^N_1 \fc{1}{s} \left| \int^\iy_1 t^{-1-r}
 \cos \left(st - \fc{\pi n}{2}\right)\ dt\ \right|\ ds
=\int^N_1 \fc{1}{s} \left| \fc{1}{s} t^{-1-r} \sin
\left(st - \fc{\pi n}{2}\right)\right|^\iy_1\\ \quad\\
&+\fc{1+r}{s} \int^\iy_1 t^{-2-r}\bigm| \sin \left(st - \fc{\pi n}{2}\right)\
dt\bigm|\ ds\leq C \int^N_1 \df{ds}{s^2},\endalign$$ \mk\flushpar
and Example 2 is proved.\endexample

\bk

\subheading{5. RADIAL FUNCTIONS WITH CONVEXITY CONDITIONS}
\sk
{\bf 5.1.} An asymptotics of the Fourier transform may be established not very
often. Some special conditions, like convexity, must be layed on a function.
G. E. Shilov (see e.g., \cite{Ba}, p. 632) was the first who studied the
Fourier coefficients of convex functions. In [T2] an asymptotics of the
Fourier transform of a convex function is given in a sharp form as follows:
\sk
\proclaim{Theorem A2}\ If $f$ is convex on $[a,b],$ where $-\iy < a < b
\leq + \iy,$ and $|f'(b)| < \iy,$ then for each $r\in \BR^1,\ |r| \geq
2,$
$$\int^b_a f(t)\ e^{-irt}\ dt = \fc{i}{r}\ \left\{ f(b)\ e^{-ibr} -
f\left(a+\fc{d}{|r|}\right)\ e^{-iav}\right\} + \th\g (|r|),$$
where $d=\min \{b-a,\pi\},$ $|\th| \leq C,$ and $\g$ is monotone decreasing
so that $$\int^\iy_2 \g(t)\ dt \leq \fc{1}{d}\ V_f + |f' (b)|.$$ \ep
\sk
{\bf 5.2.} In \cite{T3}, \cite{T4}, \cite{L1}, this theorem was generalized
to the radial case. Here we prove the stronger result formulated in \cite{BL1}.

\sk
\proclaim{Theorem 3}  Let a radial function $f$ satisfy \rom{(1)}.
Suppose further that $f_0$ is supported on $[0,1]$ and $F$ is
continuous on $[0,\infty)$ and convex on $[0,1]$. Then for $|x|=r\geq 2$
$$\hf(x)=2^{\fc{n+1}{3}}\pi^{\fc{n-1}{2}}(-1)^{\left[\fc{n}{2}\right]}r^{-n}F
\left(1-\fc{1}{r}\right)\cos\left(r-\fc{\pi n}{2}\right)+\th\g(r),$$
where $|\th|\le C$ and $\g$ decreases monotonously so that
$$\int_2^\iy r^{n-1}|\g(r)|\,dr\le V_F.$$   \ep
\pf Using (6) we obtain
$$\align \hf(x)&=\hf_0(r)=\fc{(2\pi)^{\fc{n}{2}}(-1)^{\left[\fc{n}{2}\right]}}
{\G\left(\fc{n-1}{2}\right)}r^{1-\fc{n}{2}}\int_0^1F(t)t^{\fc{n}{2}}Q(rt)dt\\
\quad\\  &=\fc{(2\pi)^{\fc{n}{2}}(-1)^{\left[\fc{n}{2}\right]}}{\G
\left(\fc{n-1}{2}\right)} r^{1-\fc{n}{2}}\left\{F\left(1-\fc{1}{r}\right)
\int_0^1t^{\fc{n}{2}}Q(rt)dt\right.\\ \quad\\
&+\int_{1-\fc{1}{r}}^1\left[F(t)-F\left(1-\fc{1}{r}\right)\right]t^{\fc{n}{2}}
Q(rt)dt\\ \quad\\
&-\fc{1}{r}\left.\int_{\fc{1}{r}}^{1-\fc{1}{r}}F'(t)t^{\fc{n}{2}}q(rt)dt
+\int_0^{\fc{1}{r}}\left[F(t)-F\left(\fc{1}{r}\right)\right]t^{\fc{n}{2}}
Q(rt)dt\right\}.\endalign$$
Lemma 1 	and (7) yield
$$\align & \int_0^1t^{\fc{n}{2}}Q(rt)dt=\fc{1}{r}q(r) \tag16\\
& =\G\left(\fc{n-1}{2}\right)r^{-\fc{
n+1}{2}}J_{n-\fc{1}{2}}(r)+\zeta_n r^{-\fc{n+2}{2}}+O\left(
r^{-\fc{n+4}{2}}\right). \endalign$$ We denoted $\z_{{n-1\over 2},n}$
by $\z_n.$ The remainder term may be obviously referred to as $\g.$
In order to obtain the main term of asymptotics apply (8) to
$J_{n-\fc{1}{2}}(r).$
Then, using Lemma 2 for the upper estimate, we obtain
$$\align &\left|\int_{1-{1\over r}}^1 [F(t)-F(1-{1\over r})]
t^{{n\over 2}} Q(rt)\,dt\right|
=\left|\int_{1-\fc{1}{r}}^1 t^{\fc{n}{2}}Q(rt)dt\int_{1-\fc{1}{r}}^t
F'(s)ds\right|\\ \quad\\
&\quad
=\left|\int_{1-\fc{1}{r}}^1F'(s)ds\int_s^1t^{\fc{n}{2}}Q(rt)dt\right|\le
Cr^{-\fc{n}{2}}\int_{1-\fc{2}{r}}^1|F'(s)|(1-s)ds.\tag17\endalign$$
Let us verify that the right-hand side satisfies conditions which
are defining $\g.$ Monotonicity is obvious. Further
$$\int_2^\iy
r^{n-1}r^{1-\fc{n}{2}}r^{-\fc{n}{2}}dr\int_{1-\fc{2}{r}}^1|F'(s)|(1-s)
ds=\int_0^1|F'(s)|(1-s)ds\int_2^{\fc{1}{1-s}}dr\le 2V_F.$$
Lemma 1 yields for large $rt$
$$t^{\fc{n}{2}}q(rt)=\G\left(\fc{n-1}{2}\right)t^{\fc{1}{2}}r^{-\fc{n-1}{2}}
 J_{n-\fc{1}{2}}(rt)+\zeta_n r^{-\fc{n}{2}}+O
 \left(\fc{1}{t}r^{-\fc{n+2}{2}} \right).$$
For $rt$ small the remainder term is $O\left(t^{\fc{n}{2}}\right).$ Then
$$\align-\fc{1}{r}&\int_{\fc{1}{r}}^{1-\fc{1}{r}}F'(t)t^{\fc{n}{2}}q(rt)dt+
\int_0^\fc{1}{r}\left[F(t)-F\left(\fc{1}{r}\right)\right]t^{\fc{n}{2}}Q(rt)
dt\\ \quad\\
&=-\G\left(\fc{n-1}{2}\right)r^{-\fc{n+1}{2}}\int_{\fc{1}{r}}^{1-\fc{1}{r}}
F'(t)t^{\fc{1}{2}}J_{n-\fc{1}{2}}(rt)dt-\z_n r^{-\fc{n+2}{2}}
F\left(1-\fc{1}{r}\right)\\
&\quad +\z_n r^{-\fc{n+2}{2}}F\left(\fc{1}{r}\right)+O\left(
r^{-\fc{n+4}{2}}\int_{\fc{1}{r}}^{1-\fc{1}{r}}|F'(t)|\fc{dt}{t}\right)
-\fc{1}{r}\int_0^{\fc{1}{r}}F'(t)t^{\fc{n}{2}}q(rt)dt.\tag18\endalign$$
But we have
$$\align -\fc{1}{r}\int_0^{\fc{1}{r}}F'(t)t^{\fc{n}{2}}q(rt)dt&=-\G\left(
\fc{n-1}{2}\right)r^{-\fc{n+1}{2}}\int_0^{\fc{1}{r}}F'(t)t^{\fc{1}{2}}J_{n -
\fc{1}{2}}(rt)dt\\ \quad\\
&-\z_n r^{-\fc{n+2}{2}}F\left(\fc{1}{r}\right)+O\left(
\fc{1}{r}\int_0^{\fc{1}{r}}|F'(t)|t^{\fc{n}{2}}dt\right).\endalign$$
Now the second summand annihilates the third summand on the right-hand side of
(18), and after applying (8) the first summand and the remainder term are
estimated by
$$r^{-\fc{n}{2}}\int_0^{\fc{2}{r}}|F'(t)|tdt.\tag19$$
Integrability of (19) over $[2,\iy]$ with the weight $r^{1-{n\over 2}}r^{n-1}$
is verified like for (17).
The second summands in (16) and (18) are mutually annihilated as well.
Let us prove that the remainder term in (18) is a part of $\g.$
$$\align &\int_2^\iy r^{n-1}dr \
r^{-n-1}\int_{\fc{1}{r}}^{1-\fc{1}{r}}|F'(t)|\fc{dt}{t}=
\int_2^\iy \fc{dr}{r^2}\left(\int_{\fc{1}{r}}^{\fc{1}{2}}+\int_{
\fc{1}{2}}^{1-\fc{1}{r}}\right)|F'(t)|\fc{dt}{t}\\ \quad\\
&\qquad \qquad \le\int_0^{\fc{1}{2}}|F'(t)|\fc{dt}{t}\int_{\fc{1}{t}}^\iy
\fc{dr}{r^2}+2\int_{\fc{1}{2}}^1|F'(t)|dt\le 2\int_0^1|F'(t)|dt.\endalign$$
It remains to estimate the first summand on the right-hand side of (18).
Integrating by parts we obtain
$$\align
\int_{\fc{1}{r}}^{1-\fc{1}{r}}F'(t)t^{\fc{1}{2}}J_{n-\fc{1}{2}}(rt) dt&=
F'(t)\left.\int_0^ts^{\fc{1}{2}}J_{n-\fc{1}{2}}(rs)ds
\right|_{\fc{1}{r}}^{1-\fc{1}{r}}\\ \quad\\
&-\int_{\fc{1}{r}}^{1-\fc{1}{r}}\biggl[\int_0^ts^{\fc{1}{2}}
J_{n-\fc{1}{2}}(rs)ds\biggr] \ dF'(t).\endalign$$
Let us estimate the integrated terms. The integral is estimated by Lemma 3.
Since $F'$ can be assumed monotone the integrated term for $t=1-\fc{1}{r}$
is estimated by (17), and for $t=\fc{1}{r}$ by (19).
Again, using Lemma 3 we have
$$\left|\int_{\fc{1}{r}}^{1-\fc{1}{r}}\int_0^ts^{\fc{1}{2}}J_{n-\fc{1}{2}}(rs)
ds\ dF'(t)\right|\le
C\fc{1}{r}\int_{\fc{1}{r}}^{\fc{1}{2}}t^{\fc{1}{2}}|dF'(t)| +
Cr^{-\fc{1}{2}}\int_{\fc{1}{2}}^{1-\fc{1}{r}}|dF'(t)|.$$
Since $F'(t)$ is monotone, then
$$\left.\int_{\fc{1}{r}}^{\fc{1}{2}}t^{\fc{1}{2}}|dF'(t)|=\left|\int_{\fc{1}{r}}^{\fc{1}{2}}
t^{\fc{1}{2}}dF'(t)\right|=\left|F'(t)t^{\fc{1}{2}}\right|_{\fc{1}{r}}^{\fc{1}{2}}
-\fc{1}{2}\int_{\fc{1}{r}}^{\fc{1}{2}}t^{-\fc{1}{2}}F'(t)dt\right|$$
and
$$\int_{\fc{1}{2}}^{1-\fc{1}{r}}|dF'(t)|=
\left|F'(1-\fc{1}{r})-F'\left(\fc{1}{2}\right)\right|.$$
Again the integrated terms are estimated by (17), (19), and also by
$$\left|F'\left(\fc{1}{2}\right)\right|\int_2^\iy
r^{n-1}r^{-n-\fc{1}{2}}dr=\sqrt 2\left|F'\left(\fc{1}{2}\right)\right|.$$
Taking into account the factor $r^{-n-\fc{1}{2}}$ we have for the integral the
final estimate
$$r^{-n-\fc{1}{2}}\int_{\fc{1}{r}}^{\fc{1}{2}}t^{-\fc{1}{2}}|F'(t)|dt.$$
It remains to verify that this value is also suitable for $\g.$
Again the monotonicity is obvious. Further
$$\align & \int_2^\iy
r^{n-1} r^{-n-\fc{1}{2}}dr\int_{\fc{1}{r}}^{\fc{1}{2}}t^{-\fc{1}{2}}|F'(t)|dt\\
\quad\\ =&\int_0^{\fc{1}{2}}|F'(t)|t^{-\fc{1}{2}}dt\int_{\fc{1}{2}}^\iy
r^{-\fc{3}{2}}\,dr = 2\int_0^{\fc{1}{2}}|F'(t)|dt\le 2V_F.\endalign$$
The theorem is completely proved. \qed\edm
\remark{Remark 4} This theorem also shows that convexity conditions
posed on $F$ allow to avoid condition (14). \endremark
\sk
{\bf 5.3.} A criterion of integrability of Fourier transforms is a simple
corollary to this theorem.
\sk
\proclaim{Corollary 2}
Under conditions of Theorem \rom{3} the Fourier transform
$\hf$ is integrable if and only if the integral
$$\int_0^1\fc{f_0(t)}{(1-t)^{\fc{n+1}{2}}}dt\tag20$$ converges.\ep
\sk

\pf It is obvious that $\hf\in L^1(\Bbb R^n)$ is equivalent to the finiteness
of
$$\int_0^\iy r^{n-1}r^{-n}\left|F\left(1-\fc{1}{r}\right)\right|dr =
\int_0^\fc{1}{2}|F(1-t)|\fc{dt}{t}.$$
Recalling the definition of $F$ and taking into account that by convexity
$F(1-t)$ preserves a sign near the origin, we get the following integral to
estimate:
$$\int_0^{\fc{1}{2}}\fc{1}{t}f_0^{\left(\fc{n-1}{2}\right)}(1-t)dt.$$
If $n$ is odd the usual integration by parts leads to the required result.
When $n$ is even it suffices to consider
$$\int_0^{\fc{1}{2}}t^{-\fc{n}{2}}f_0^{\left(\fc{1}{2}\right)}(1-t)dt=
\int_0^\iy t^{-\fc{n}{2}}\fc{d}{dt}\int_{1-t}^1(s-1+t)^{-\fc{1}{2}} f_0(s)ds \
dt.$$
Integrate by parts in $t.$ We get
$$\align \int_0^{\fc{1}{2}} t^{-\fc{n}{2}} f_0^{\left(\fc{1}{2}\right)}
(1-t)dt & =t^{-\fc{n}{2}}\left.
\int_{1-t}^1 (s-1+t)^{-\fc{1}{2}}f_0(s)ds\right|_0^{\fc{1}{2}} \\ \quad\\
&+\fc{n}{2}\int_0^{\fc{1}{2}}t^{-1-\fc{n}{2}}\int_{1-t}^1 (s-1+t)^{-\fc{1}{2}}
f_0(s)ds\ dt. \endalign$$
The integrated terms are bounded because the integral (20) converges.
After changing the order of integration, the last integral is equal to
$$\align &\int_{\fc{1}{2}}^1f_0(s)ds\int_{1-s}^{\fc{1}{2}}(s-1+t)^{-\fc{1}{2}}
t^{-1-\fc{n}{2}}dt\\ \quad\\
=&\int_{\fc{1}{2}}^1f_0(s)ds\int_{1-s}^\iy(s-1+t)^{-\fc{1}{2}}t^{-1-\fc{n}{2}}\
dt-\int_{\fc{1}{2}}^1t_0(s)ds\int_{\fc{1}{2}}^\iy (s-1+t)^{-\fc{1}{2}}
t^{-1-\fc{n}{2}}dt.\endalign$$
The last integral is obviously bounded.
An in the first one the inner integral is the Weyl integral of order
$\df{1}{2}$ of the function $t^{-1-\fc{n}{2}}.$
By the formula in \cite{BE2}, p.201 it is equal to
$\fc{\G\left(\fc{n+1}{2}\right)}{\G\left(\fc{n+2}{2}\right)}
(1-t)^{-\fc{n+1}{2}},$
and the corollary is established. \qed\edm

\subheading{6. ACKNOWLEDGEMENTS}

The main problems of this work were posed to the author by R.M. Trigub.
The author thanks him sincerely for this and for
many helpful discussions. Also the author is grateful to E. Belinsky and
W. Trebels for many suggestions which made the presentation better.
\bk

\enddocument